\documentclass[aps,pre,reprint,groupedaddress]{revtex4-1}
\usepackage{mathtools}
\usepackage{amssymb}
\usepackage{amsmath, amsthm, amssymb}
\usepackage{graphicx}
\usepackage{float}
\usepackage{color}
\usepackage{mathrsfs}
\usepackage{enumitem}
\usepackage[ngerman,english]{babel} 
\newcommand{\R}{\mathbb{R}}
\newcommand{\Rn}{\mathbb{R}^n}

\newcommand{\e}{\text{e}}
\newcommand{\C}{\mathbb{C}}
\newcommand{\Cn}{\mathbb{C}^n}
\newcommand{\dhex}{$H^-$}
\newcommand{\uhex}{$H^+$}
\newcommand{\s}{$S$}
\newcommand{\ep}{\varepsilon}
\providecommand{\abs}[1]{\left\vert#1\right\vert}

\providecommand{\ali}[1]{\begin{align}#1\end{align}}
\providecommand{\alis}[1]{\begin{align}\begin{split}#1\end{split}\end{align}}
\providecommand{\alinon}[1]{\begin{align*}#1\end{align*}}
\providecommand{\figref}[1]{Fig.\ref{#1}}
\providecommand{\figrefa}[1]{Fig.\ref{#1}(a)}
\providecommand{\figrefb}[1]{Fig.\ref{#1}(b)}
\providecommand{\figrefc}[1]{Fig.\ref{#1}(c)}

\providecommand{\figrefj}[1]{Fig.\ref{#1}(j)}

\date{\today}

\begin{document}
\title{Tristability between stripes, up-, and down-hexagons and snaking bifurcation branches of spatial connections between up- and down-hexagons }
\author{D. Wetzel}
\affiliation{Institut f\"ur Mathematik,\\ Carl-von-Ossietzky Universit\"at Oldenburg,\\ 26111 Oldenburg, Germany.}
\begin{abstract}
Third order amplitude equations on hexagonal lattices can be used for predicting the existence and stability of stripes, up- and down-hexagons in pattern forming systems. These amplitude equations predict the nonexistence of bistable ranges between up- and down-hexagons and tristable ranges between stripes, up- and down-hexagons. In the present work we use fifth order amplitude equations for finding such bistable and tristable ranges for a generalized Swift-Hohenberg equation and discuss stationary front connections between up- and down-hexagons.
\end{abstract}
\maketitle
\section{Introduction}
It was shown by Turing \cite{turing} in 1952 that nonhomogeneous steady states arise in reaction-diffusion systems when a homogeneous state is unstable for the full system and stable for the kinetics. This discovery was followed by a large number of works, where systems of different scientific disciplines such as biology \cite{gierermeinhardt,murray}, chemistry \cite{PP00,epstein03,epstein06,K06,epstein09,epstein3d11}, ecology \cite{meron01,meron04,meron05,MS06}, and physics \cite{swift77,goswkn84} are studied for so-called Turing patterns. 

Typical 2D Turing patterns are labyrinth patterns, gaps, and spots. If spots and gaps have a hexagonal structure, they are referred to as up- and down-hexagons, respectively. A special subset of labyrinth patterns are stripes. Stability transitions for such patterns are studied via third order amplitude equations in \cite{TransPattern}. These amplitude equations predict that up- and down-hexagons corresponding to the same wavenumber cannot be stable at the same time. In contrast, tristable ranges between stripes, up-, and down-hexagons for a generalized Swift-Hohenberg model and for a specific system modeling vegetations of semiarid ecosystems are found in \cite{hilaly} and \cite{meron04} via using numerical time integration, respectively.
That tristable ranges between stripes, up-, and down-hexagons can be predicted using fifth order amplitude equations is shown in \cite{hilaly}. 

During the last 30 years, a great interest arose in localized Turing patterns. 
It was already understood by Pomeau \cite{pomeau} in 1986 that localized patterns on homogeneous backgrounds can be found in reaction-diffusion systems, when both corresponding states are stable. If the  steady system has a spatial conserved quantity, it must necessarily be equal for both states  in order to permit the existence of such a connecting state. The point where this takes place is called the Maxwell point. Via changing the wavelength of the patterned state, one can find more Maxwell points. As a consequence, branches of such states are found which move back and forth in parameter space and pass stable and unstable ranges. This scenario is referred to as snaking \cite{woods}. There are a lot of works, which investigate this effect over 1D domains (see e.g. \cite{burke,bukno2007,BKLS09}). For a detailed analysis  using the Ginzburg-Landau formalism and beyond all order asymptotics see \cite{chapk09,dean11}.

Localized patterns correspond to two stationary fronts, which are glued together. Thus, such patterns exist on unbounded domains so that one cannot find these states on bounded domains. What remains are stationary states, which are periodic in space and for which the corresponding orbits pass near the homogeneous state and the Turing pattern. We also call such states localized or connecting patterns.
Their branches also show a snaking behavior (see \cite{bbkm2008,dawes08,dawes09,KAC09,hokno2009} for further details).

It was also understood by Pomeau \cite{pomeau} that stationary connections between hexagons and stripes should exist in ranges where both patterns are stable.
The following connecting patterns are observed in \cite{hilaly} using numerical time integrations:
\begin{enumerate}[label=(\alph*)]
\item localized patches of hexagons on a homogeneous background,
\item connections between hexagons and stripes,
\item connections between up- and down-hexagons.
\end{enumerate}

Auto \cite{auto97b} is a software, which exists for a long time and which is often used for calculating solution branches of 1D patterns numerically and especially snaking branches. For a long time there was no software for performing these calculations for 2D patterns, which explains why numerically calculated branches of 2D patterns are rare. In \cite{hexsnake} a custom continuation code was used for calculating snaking branches corresponding to states of type (a). Since a couple of years the software pde2path \cite{p2phome,p2pure} exists, which is designed for calculating solution branches of PDEs over 1D, 2D, and 3D domains. This software is used in \cite{schnaki,w16} for calculating branches corresponding to states of type (a) and (b). In the present work we also use this software for investigating branches of solutions of type (c), which has not been reported before.

The present work is organized as follows. We start with a quick review of third order amplitude equation reductions. After this we show a bifurcation diagram for which we used numerical continuation methods for following the stripe, gap, and spot branches bifurcating from a homogeneous state of the vegetation model mentioned above and observe a tristable range as well. In order to predict tristable ranges via the amplitude formalism we use a generalized Swift-Hohenberg equation, which is scaled with a parameter $\ep$ such that we are able to reduce the full equation to a system of fifth order amplitude equations. In the following a comparison of these predictions with numerical solutions is performed which shows that these amplitude equations give acceptable predictions if $\ep$ is small. At the end of the present work we show a snaking branch corresponding to connections between up- and down-hexagons and explain this behavior using conserved quantities.   

\section{Third order amplitude equations\label{o3}}
Consider
\ali{u_t=D\Delta u + f(u,\lambda),\label{rds}}  
where $u:\R \times \R \times \R_{\ge 0} \rightarrow \Rn$, $(x,y,t)\rightarrow u(x,y,t)$, $\lambda\in \R$ a control parameter, $D\in\R^{n\times n}$, $\Delta=\partial^2_x+\partial^2_y$, $x,y$ space coordinates, and $t$ the time.
Let $u^*$ be a homogeneous steady state of \eqref{rds} and $(u^*_c,\lambda_c)$ a Turing point with corresponding critical wavenumber $k_c$.  Using the ansatz
\ali{
u-u^*=(A_1 e_1+A_2 e_2+A_3 e_3)\Phi+\text{c.c.}+\text{h.o.t.} \label{ansatz}
}
with $u^*\in\Rn$,  $\Phi\in\Cn$, $e_j=\e^{i (x,y)\cdot \mathbf{k_j}}$, $A_j=A_j(t)\in\C$ for $j=1,2,3,$
$\mathbf{k_1}=k_c \ (1,0)^T$, $\mathbf{k_2}=0.5 \ k_c \ (-1,\sqrt{3})^T$, and $\mathbf{k_3}=0.5 \ k_c \ (-1,-\sqrt{3})^T$, we can 
reduce the full system to a system of amplitude equations given by
 \begin{align}
 \label{gl}
 \begin{split}
\dot{A_1}=c_1 A_1+c_2 \overline{A_2}\overline{A_3}&+c_3 A_1 |A_1|^2\\  &+c_4 A_1(|A_2|^2 +|A_3|^2),\\
\dot{A_2}=c_1 A_2+c_2 \overline{A_1}\overline{A_3} &+c_3 A_2 |A_2|^2\\ &+c_4 A_2(|A_1|^2 +|A_3|^2),\\
\dot{A_3}=c_1 A_3+c_2 \overline{A_1}\overline{A_2} &+c_3 A_3 |A_3|^2\\ &+c_4 A_3(|A_1|^2 +|A_2|^2) 
\end{split}
\end{align}
(see \cite{goswkn84,Hoyle} for further details).
By the center manifold reduction (see chapter 13 of \cite{SU17}) an $\ep_0>0$ exists such that for all $0<\ep<\ep_0$ and
\alinon{
u_A= \ &2\ep \bigg[  A_1 \cos(x)+A_2 \cos\left(-0.5\left(x+\sqrt{3}y\right)\right)\\
&+A_3 \cos\left(-0.5\left(x-\sqrt{3}y\right)\right) \bigg]
} 
with $(A_1,A_2,A_3)$ solving \eqref{gl} a solution $u_O$ of \eqref{sh} exists such that $u_A-u_O=\mathcal{O}(\ep^2)$. 

Such approximations are only valid near $\lambda_c$ and for small amplitudes. Typical time independent solutions of \eqref{gl} are up-hexagons $H^+$, stripes $S$, and down-hexagons $H^-$ for which the triple $(A_1,A_2,A_3)$ is given by
 $(h^+,h^+,h^+)$, $(s,0,0)$, $(h^-,h^-,h^-)$, respectively. Here $h^+\in \R_{>0}$, $s\in \R$, and $h^-\in \R_{<0}$. Their stability can be obtained from the Jacobian of \eqref{gl}. Unfortunately the amplitude reduction is only valid for small amplitudes (near $\lambda_c$), and if the quadratic terms of the Taylor expansion of $f$ in \eqref{rds} around $u^*$ are small. 


Stability transitions between \uhex, \s, and \dhex\ via varying $c_1$ and $c_2$ are discussed in \cite{TransPattern} and the following result is pointed out. For fixed parameters $c_1, \ c_2, \ c_3, \ c_4$ for which the three states \uhex, \s, and \dhex\ exist, the system \eqref{gl} can only predict one of the following situations:

\begin{itemize}
\item only one of the states \uhex, \s, \dhex\ is stable (mono- stable case),
\item \s\ and \uhex\ or \dhex\ are stable (bistable case),
\item all three states are unstable.
\end{itemize} 
Thus, bistable ranges between \uhex\ and \dhex\ and tristable ranges between \uhex, \s, and \dhex\ are not predictable using \eqref{gl}.

\begin{figure}
\begin{center}
\begin{minipage}{0.23\textwidth}
(a) $\beta=3$\\
\includegraphics[width=1\textwidth]{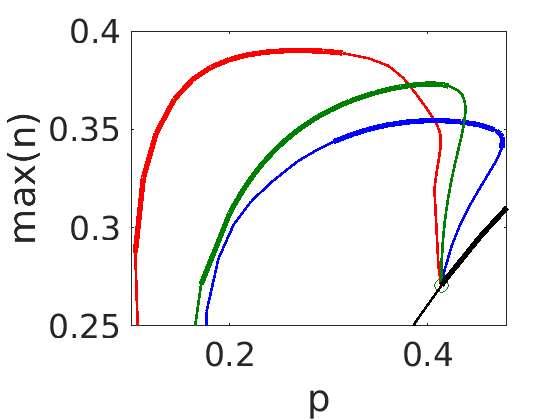}
\end{minipage}
\begin{minipage}{0.23\textwidth}
(b) $\beta=3.5$\\
\includegraphics[width=1\textwidth]{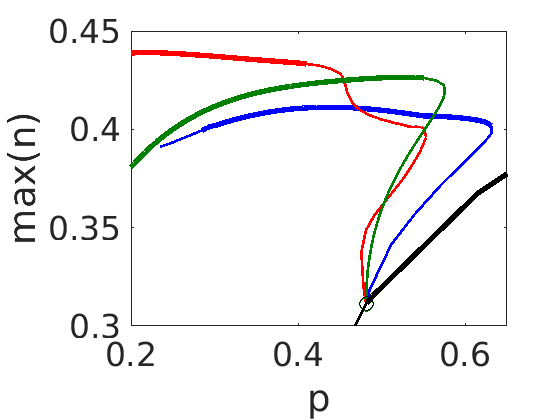}
\end{minipage}
\begin{minipage}{0.15\textwidth}
\footnotesize{(c) Up-hexagons}\\
\includegraphics[width=1\textwidth]{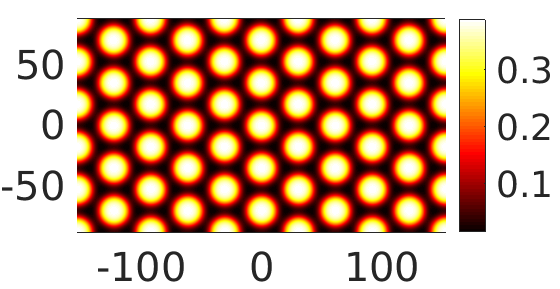}
\end{minipage} 
\begin{minipage}{0.15\textwidth}
\footnotesize{(d) Stripes}\\
\includegraphics[width=1\textwidth]{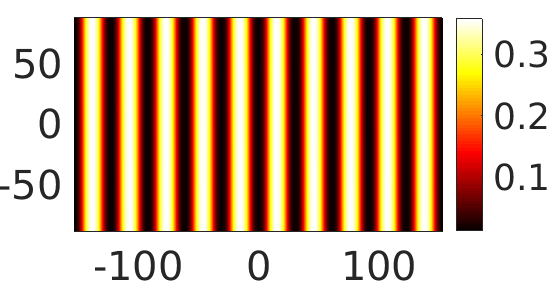}
\end{minipage} 
\begin{minipage}{0.15\textwidth}
\footnotesize{(e) Down-hexagons}\\
\includegraphics[width=1\textwidth]{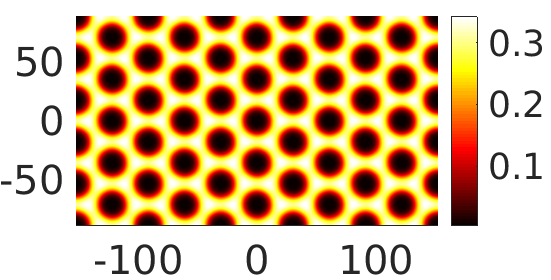}
\end{minipage}
\end{center}  
\caption{The system \eqref{meron} for the parameter set \eqref{para} has one trivial homogeneous state $(n,w)=(0,p)$, which is stable for $p<p_0$ and unstable for $p>p_0$ with $p_0\approx 0.157$. At this transition of stability a nontrivial homogeneous state bifurcates, which is stable at the beginning, becomes Turing unstable a bit later (in $p\approx 0.169$) and becomes stable again later on (in $p\approx 0.413$). The black line in (a) represents the nontrivial homogeneous states. Here and in the following thick and thin lines represent stable and unstable states, respectively. We use the continuation and bifurcation software tool pde2path \cite{p2phome,p2pure} to follow the solution branches which bifurcate in the second Turing point. Here we use the domain $\Omega=(-l_x,l_x)\times (-l_y,l_y)$ with $l_x=20\pi/k_c$, $l_y=20\pi / (\sqrt{3}k_c)$ and Neumann boundary conditions, while $k_c=0.206$ is the critical wavenumber corresponding to the second Turing point. The red, green, and blue branches in (a) and (b) correspond to $H^+$, $S$, and $H^-$, respectively. We can see in (a) that there is a small tristable range. All three states are stable for $p=0.301$ and their density plots for $n$ are shown over the domain $\Omega=(-l_x/2,l_x/2)\times (-l_y/2,l_y/2)$ in (c)-(e). We set $\beta=3.5$ for developing (b). Here we have a larger tristable range.}\label{meronplot}
\end{figure}

The following vegetation model for semi-arid ecosystems
\begin{equation}\label{meron}
\begin{aligned}
n_t&=\Delta n +\left(\frac{\gamma w}{1+\sigma w}-\nu\right)n-n^2,\\
w_t&=\delta \Delta(w-\beta n)+p-(1-\rho n)w-w^2n
\end{aligned}
\end{equation}
is discussed in \cite{meron04} and used in \cite{TransPattern} to apply the amplitude reduction and compare it with numerical results. Here $n$, $w$, and $p$ represent the vegetation density, ground water density, and precipitation, respectively. $p$ is used as a control parameter. The other parameters are given by 
\ali{\gamma=\sigma=1.6, \ \nu=0.2, \ \rho=1.5, \ \beta=3, \ \delta=100.\label{para}} 
Please see \cite{meron04} for modeling details and the meaning of these parameters.
Numerical time integration methods are used in \cite{meron04} to follow the stable parts of stripes and up- and down-hexagons by varying $p$ in both directions and to see where transitions between these patterns occur. Here, a very small tristable range between \uhex, \s, and \dhex\ is observed. In the present work we use numerical path following methods to calculate the corresponding branches and observe this tristable range as well (see \figrefa{meronplot}). Furthermore, we use $\beta=3.5$ instead of $\beta=3$ and see that the tristable range is larger for this parameter set (see \figrefb{meronplot}).

We also see in \figref{meronplot} that the stripes bifurcate subcritically, turn around in a fold, and become stable after this fold. Since the amplitude equations \eqref{gl} are only expanded up to third order terms, the fold and consequently the stable stripes cannot be described via \eqref{gl}. The parameter $\delta$ is changed in \cite{TransPattern} such that both Turing points come closer together and stripes bifurcate supercritically, which makes the amplitude reduction valid for the whole existence regions of \uhex, \s, and \dhex.

\section{Fifth order amplitude equations \label{o5}}
We have seen above that the third order amplitude equation reduction on a hexagonal lattice does not predict the tristability of \uhex, \s, and \dhex, but a tristable range for the reaction-diffusion system \eqref{meron} is found numerically on bounded domains. In the present section we show that tristable ranges can be predicted via fifth order amplitude equations, which is also done in \cite{hilaly}. As for the third order case such an expansion is only valid for special systems for which the quadratic terms of the Taylor expansion of the kinetics around the homogeneous solution are small. Unfortunately this is not the case for \eqref{meron} with the corresponding parameter set \eqref{para} so that we consider the following generalized Swift-Hohenberg equation
\alis{u_t= & \ [\ep^4 c_1-(1+\Delta)^2]u\\
&+\ep^3 c_2u^2+\ep^2 c_3u^3+\ep c_4 u^4 + c_5 u^5\label{sh}}
with $c_1$, $c_2$, $c_3$, $c_4$, $c_5\in\R$ and $\ep\in \R_{>0}$. 


\eqref{sh} has the trivial solution $u=0$ which has a Turing bifurcation in $c_1=0$. Furthermore, \eqref{sh} is scaled with $\ep$ such that a fifth order amplitude equation reduction is possible and valid for small $\ep$. To do so we use the ansatz \eqref{ansatz} with $\Phi=\ep$ and $u^*=0$ and end up with the following amplitude equation system (derived at $\ep^5$)

\alis{\label{gl5}
\dot{A_1}=c_1 A_1&+2c_2 \overline{A_2}\overline{A_3} +3c_3 A_1 |A_1|^2\\ &+6c_3 A_1(|A_2|^2 +|A_3|^2)+D_1,\\
\dot{A_2}=c_1 A_2&+2c_2 \overline{A_1}\overline{A_3} +3c_3 A_2 |A_2|^2\\ &+6c_3 A_2(|A_1|^2 +|A_3|^2)+D_2, \\
\dot{A_3}=c_1 A_3&+2c_2 \overline{A_1}\overline{A_2} +3c_3 A_3 |A_3|^2\\ &+6c_3 A_3(|A_1|^2 +|A_2|^2)+D_3
}
with $D_i=12 \ c_4 \ C_{4i}+ c_5 \ C_{5i}$ for $i=1,2,3$,
\alinon{
C_{41}&=A_1^2A_2A_3+(2\abs{A_1}^2+\abs{A_2}^2+\abs{A_3}^2)\overline{A_2A_3},\\
C_{42}&=A_2^2A_1A_3+(\abs{A_1}^2+2\abs{A_2}^2+\abs{A_3}^2)\overline{A_1A_3},\\
C_{43}&=A_3^2A_1A_2+(\abs{A_1}^2+\abs{A_2}^2+2\abs{A_3}^2)\overline{A_1A_2},
}
and 
\alinon{
C_{51}= \ &30\overline{A}_1 \overline{A}_2^2 \overline{A}_3^2+10A_1\big(\abs{A_1}^4+3\abs{A_2}^4+3\abs{A_3}^4\\
&+6\abs{A_1}^2\abs{A_2}^2+6\abs{A_1}^2\abs{A_3}^2+12\abs{A_2}^2\abs{A_3}^2\big),\\
C_{52}= \ &30\overline{A}_2\overline{A}_1^2\overline{A}_3^2+10A_2\big(3\abs{A_1}^4+\abs{A_2}^4+3\abs{A_3}^4\\
&+6\abs{A_1}^2\abs{A_2}^2+12\abs{A_1}^2\abs{A_3}^2+6\abs{A_2}^2\abs{A_3}^2\big),\\
C_{53}= \ &30\overline{A}_3\overline{A}_1^2\overline{A}_2^2+10A_3\big(3\abs{A_1}^4+3\abs{A_2}^4+\abs{A_3}^4\\
&+12\abs{A_1}^2\abs{A_2}^2+6\abs{A_1}^2\abs{A_3}^2+6\abs{A_2}^2\abs{A_3}^2\big).
}

\begin{figure}
\begin{center}
\begin{minipage}{0.35\textwidth}
(a)\\
\includegraphics[width=1\textwidth]{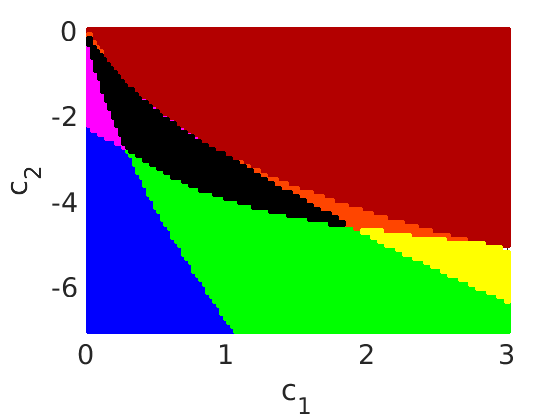}
\end{minipage}
\begin{minipage}{0.35\textwidth}
(b) $\ep=0.5$\\
\includegraphics[width=1\textwidth]{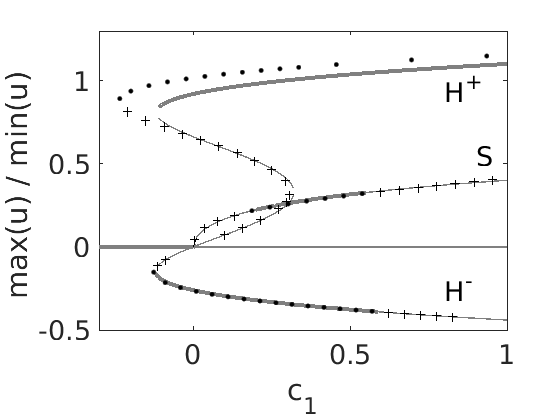}
\end{minipage}
\begin{minipage}{0.35\textwidth}
(c) $\ep=1$\\
\includegraphics[width=1\textwidth]{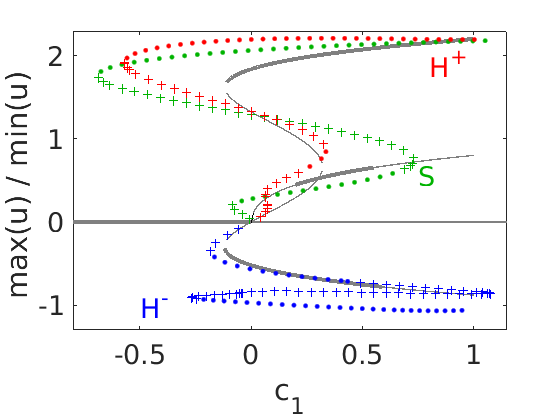}
\end{minipage}
\end{center}
\caption{We set $c_3=-1$, $c_4=5$, $c_5=-2$ for creating the plots in the present figure. We vary $c_1$ and $c_2$ and determine \uhex, \s, and \dhex\ and their stability via \eqref{gl5} for creating (a). \uhex, \s, and \dhex are monostable in the red, yellow, and blue regions, respectively. We have bistable ranges between \dhex\ and \s, \dhex and \uhex, \uhex\ and \s\ in the green, violet, and orange regions, respectively. In the black region all three states are stable. We set $c_2=-2$ for creating (b) and (c) and use $c_1$ as a control parameter. Here we use the amplitude equations \eqref{gl5} for creating the solid lines. $\bullet$ and + represent stable and unstable states, respectively, which we found numerically via the finite element and Newton's method on the domain $\Omega=(-l_x,l_x)\times (-l_y,l_y)$ with $l_x=10\pi$, $l_y=10\pi / \sqrt{3}$ and Neumann boundary conditions. Here we use $\ep=0.5$ and $\ep=1$ for (b) and (c), respectively. The vertical axis shows the maximum of $u$ for stripes and up-hexagons and the minimum of $u$ for down-hexagons.}
\label{stabregion}
\end{figure}

To see if \uhex, \s, and \dhex\ are stable or not in their existence range, we split the amplitudes in \eqref{gl5} into real and imaginary parts and gain a system with 6 equations and thus a steady state of \eqref{gl5} has 6 eigenvalues. We fix $c_3, \ c_4, \ c_5$ and determine regions of stability depending on $c_1$ and $c_2$ (see \figrefa{stabregion}). Besides monostable and bistable ranges between stripes and hexagons, we found tristable ranges between \uhex, \s, and \dhex and bistable ranges between \uhex\ and \dhex. We also fix $c_2$ and use pde2path for comparing analytical and numerical solutions, while $c_1$ is used as a control parameter (see \figrefb{stabregion} and (c)).
 
Here we see that \eqref{gl5} gives reasonable predictions for $\ep=0.5$. This does not hold for $\ep=1$. Here the stripe branch bifurcates supercritically and has a fold in $c_1\approx 10^{-6}$, where it turns around. This cannot be seen in \figrefc{stabregion}. We have some more folds on the branches corresponding to \uhex, \dhex, and \s, which cannot be predicted from \eqref{gl5}. This shows that $\ep=1$ is too large for predicting hexagons and stripes well far away and near onset, but this example gives us other interesting results. Here we have ranges where two different stable types of the same pattern type exist (a small and large amplitude pattern). This holds for all three pattern types \uhex, \dhex, and \s . Thus, we find a range where the small and large amplitude down-hexagons, the zero-solution, the small and large amplitude stripes, and the large amplitude up-hexagons are multi-stable. And furthermore another range where the small and large amplitude down-hexagons, up-hexagons, and stripes are stable.

\section{Spatial connections between patterns}
In the present section we discuss spatial connections between up- and down-hexagons solving \eqref{sh} with the parameter set used for creating \figrefb{stabregion}. We describe the behavior and the structure of the bifurcation diagram of such states, while we use the Ginzburg-Landau energy for explanations.  

\subsection{Description}
\begin{figure}
\begin{center}
\begin{minipage}{0.35\textwidth}
(a)\\
\includegraphics[width=1\textwidth]{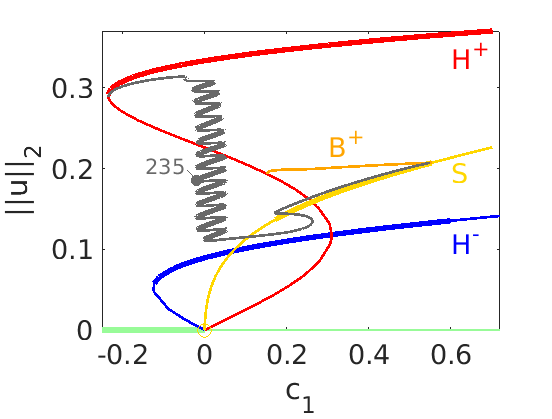}
\end{minipage}
\begin{minipage}{0.5\textwidth}
\vspace{0.3cm}
(b) Solution on $B^+$\\
\includegraphics[width=1\textwidth]{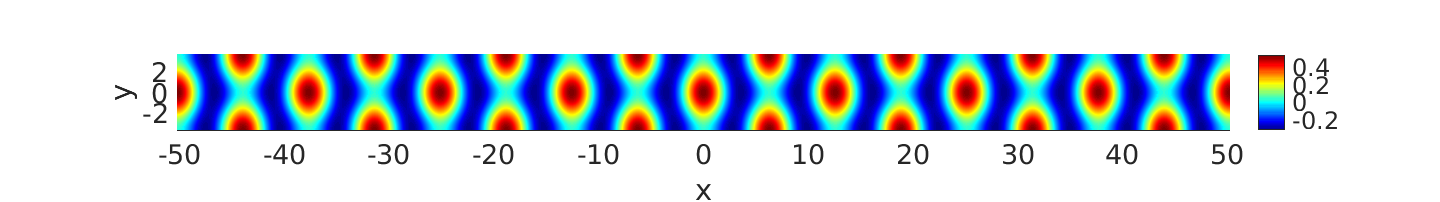}
\end{minipage}
\begin{minipage}{0.5\textwidth}
\vspace{0.3cm}
(c) Solution 235\\
\includegraphics[width=1\textwidth]{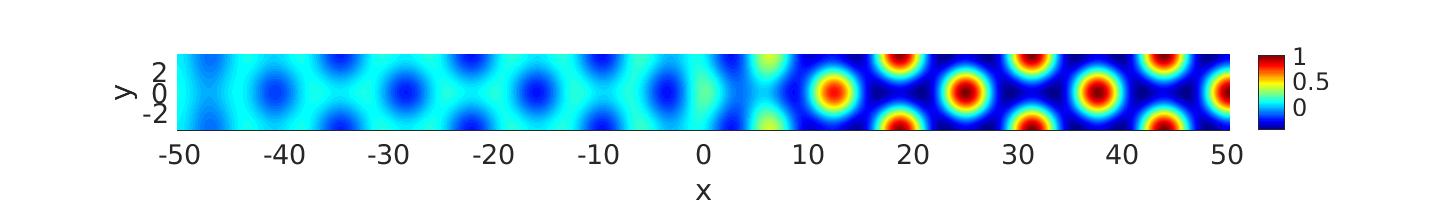}
\end{minipage}
\begin{minipage}{0.5\textwidth}
\vspace{0.3cm}
(d) left and right part of solution 235\\
\includegraphics[width=0.48\textwidth]{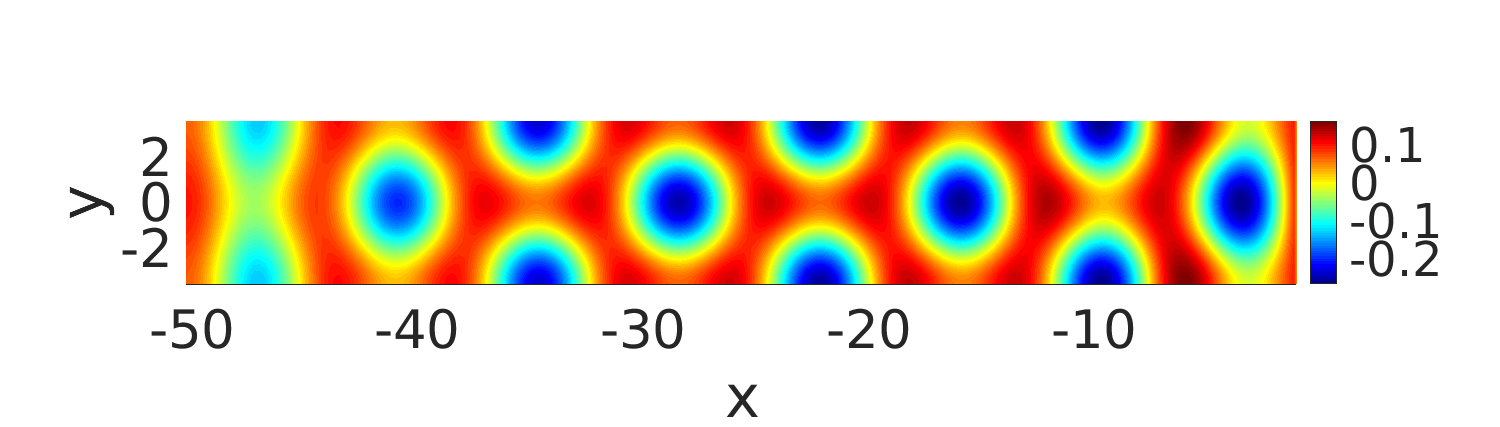}
\includegraphics[width=0.48\textwidth]{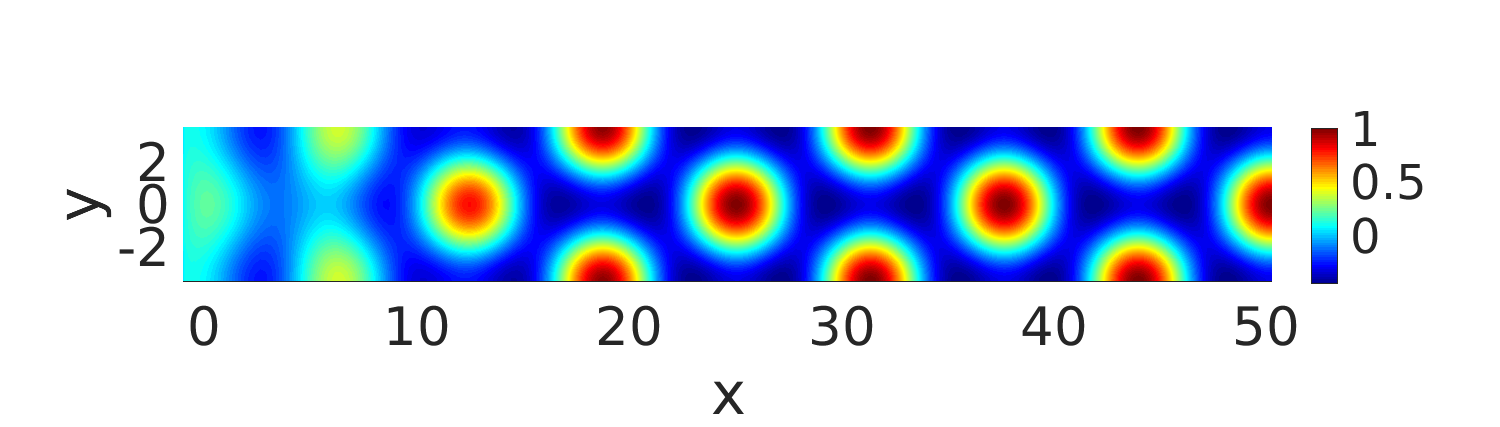}
\end{minipage}
\end{center}
\caption{(a) Numerically calculated solution branches of \eqref{sh} for $\ep=0.5$ and $c_2=-2$, $c_3=-1$, $c_4=5$, $c_5=-2$ (as for \figrefb{stabregion}). Here and in the following it holds $\|u\|_2=\sqrt{\frac{1}{\Omega}\int_{\Omega}|u|^2 \text{d}z}$, where $\Omega$ is the considered domain and $z$ represents the spatial coordinates. Here $\Omega=(-l_x,l_x)\times (-l_y,l_y)$ with $l_x=16\pi$, $l_y=2\pi / \sqrt{3}$ and we use Neumann boundary conditions. The orange branch $B^+$ is a branch of mixed mode patterns, which we call beans. This branch connects the stripe and gap branch and bifurcates in the point, where the stripes lose their stability. There is a bifurcation point on the bean branch from which the gray branch bifurcates. (b) Density plot of a solution, which lies more or less in the middle of the bean branch. (c) Density plot of a solution, which lies on the gray branch and which is labeled with 235. (d) Left and right part of solution 235 plotted seperately.} \label{hexsnake}
\end{figure}

We use a domain, which is large in the $x$-direction and small in the $y$-direction, for finding stationary solutions connecting up- and down-hexagons in space. We calculate the branches corresponding to  stripes, up- and down-hexagons on this domain and follow the branch which bifurcates in the point where the stripes lose their stability (see \figref{hexsnake}). This branch is a mixed mode pattern connecting the stripe and up-hexagons branches. We call such states beans $B^+$. They can also be found in \eqref{gl5} and are of the form $(A_1,A_2,A_3)=(A,B,B)$ with $A,B\in \R$ and $\abs{A}>\abs{B}$. A density plot of a bean solution can be found in \figrefb{hexsnake}. 

There is a bifurcation point on the bean branch located near the stripe branch from which a snaking branch bifurcates (the gray branch of \figrefa{hexsnake}). The 235th solution on this branch is labeled in \figrefa{hexsnake} (we use the labels produced by pde2path) and shown in \figrefc{hexsnake}. Here we can see that such a solution seems to be an orbit passing near up- and down-hexagons. It can be seen in \figrefc{hexsnake} and (d) that such a state changes its gap shape near the left boundary. This is because we use Neumann boundary conditions and a horizontal domain length of $n\pi$ with $n=32$ ($n$ is even).

\begin{figure}
\begin{minipage}{0.5\textwidth}
(a)\\
\includegraphics[width=1\textwidth]{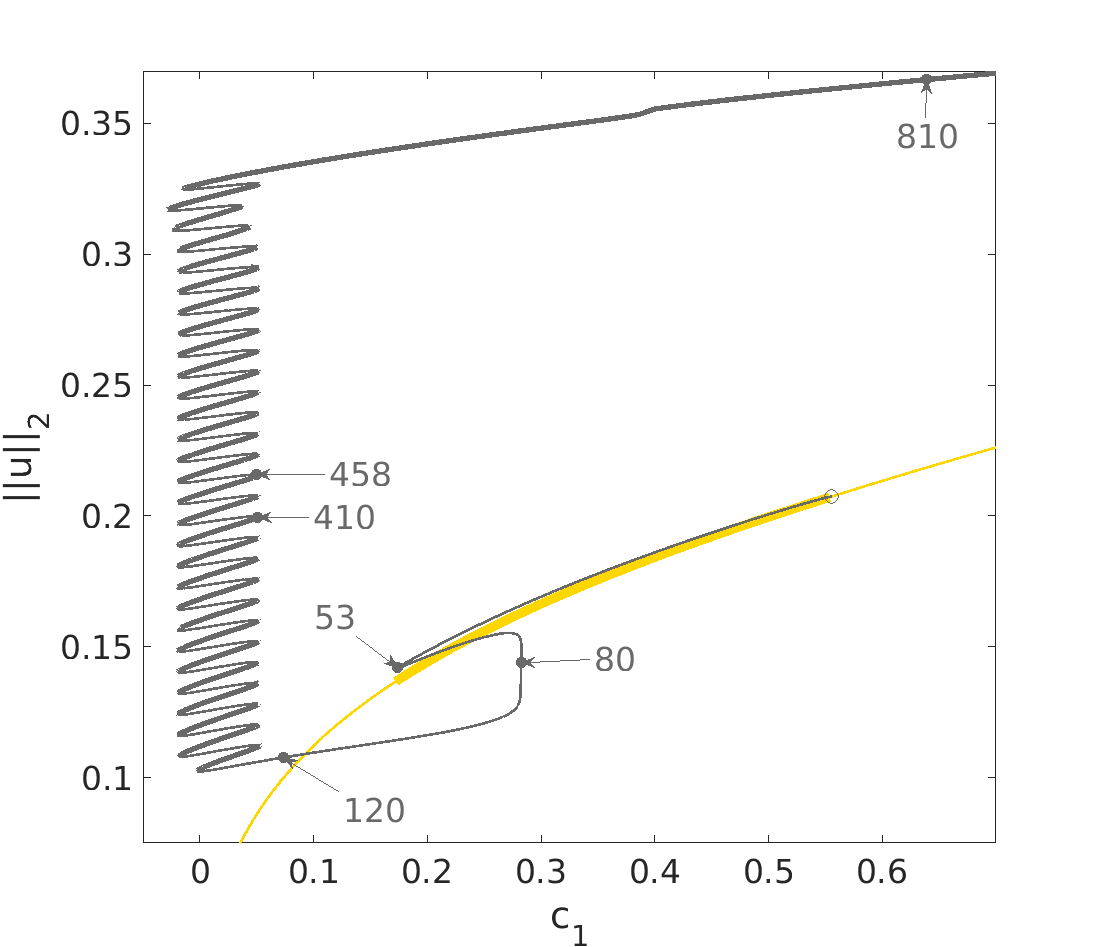}

\end{minipage}
\begin{minipage}{0.5\textwidth}
(b) Stripe solution \\
\includegraphics[width=1\textwidth]{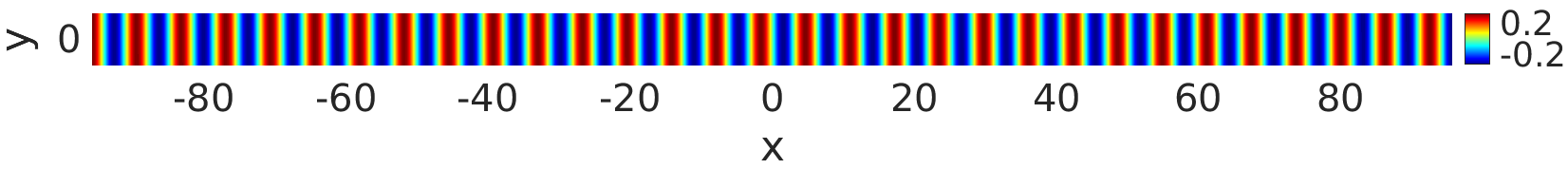}
\end{minipage}
\begin{minipage}{0.5\textwidth}
(c) 53\\
\includegraphics[width=1\textwidth]{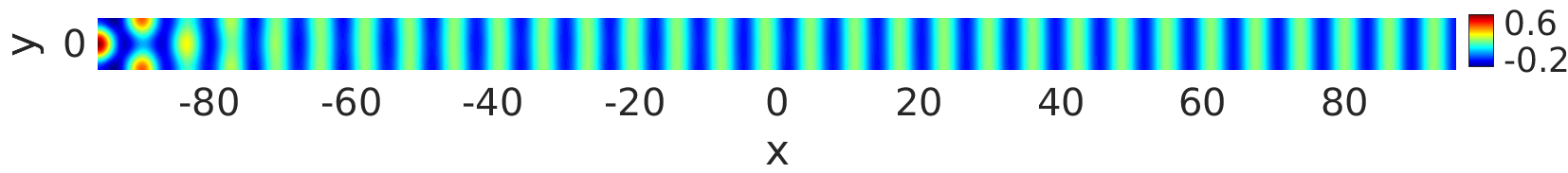}
\end{minipage}
\begin{minipage}{0.5\textwidth}
(e) 80\\
\includegraphics[width=1\textwidth]{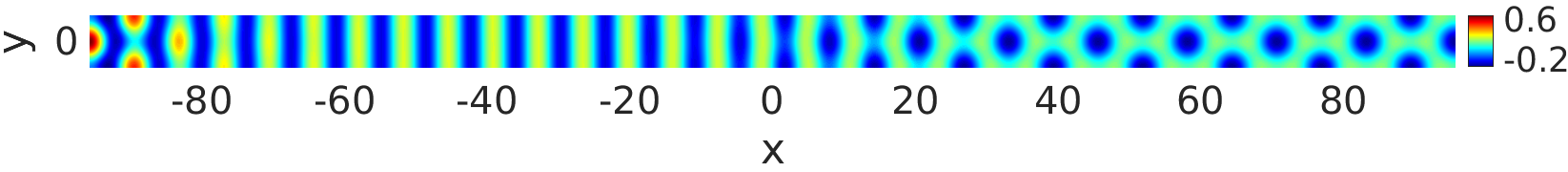}
\end{minipage}
\begin{minipage}{0.5\textwidth}
(f) 120\\
\includegraphics[width=1\textwidth]{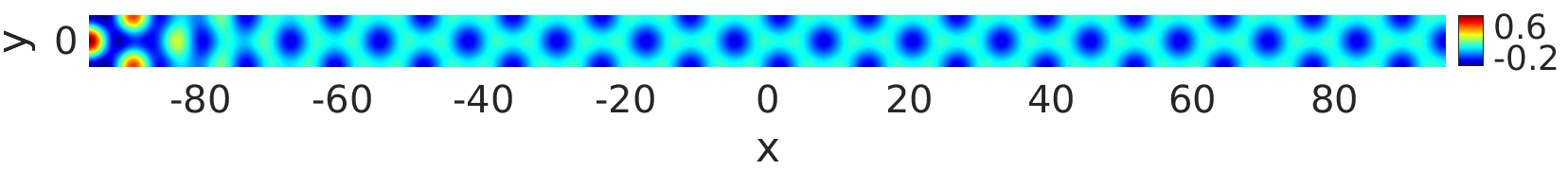}
\end{minipage}
\begin{minipage}{0.5\textwidth}
(g) 410\\
\includegraphics[width=1\textwidth]{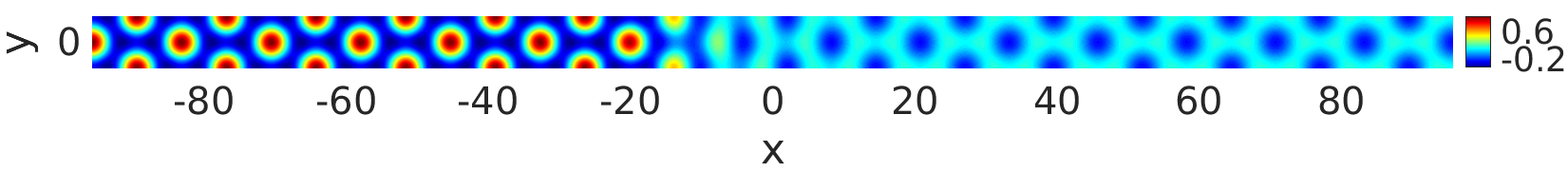}
\end{minipage}
\begin{minipage}{0.5\textwidth}
(h) 458\\
\includegraphics[width=1\textwidth]{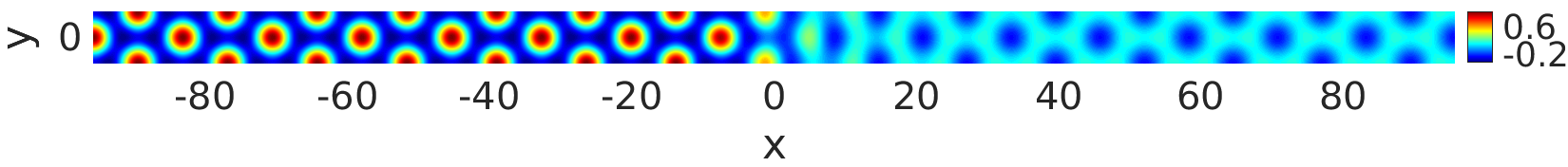}
\end{minipage}
\begin{minipage}{0.5\textwidth}
(i) 810\\
\includegraphics[width=1\textwidth]{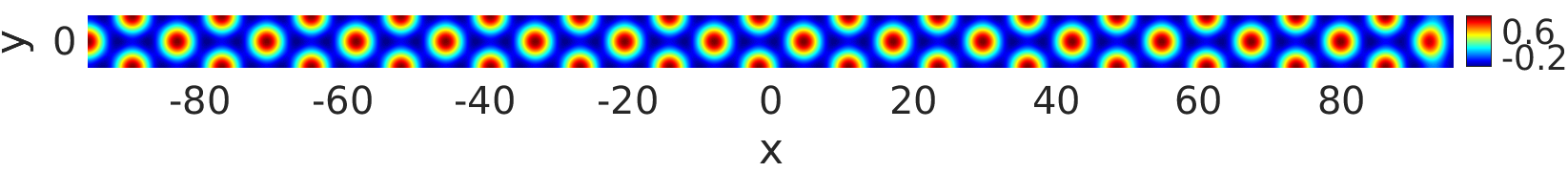}
\end{minipage}
\begin{minipage}{0.5\textwidth}
(j) left and right part of 458\\
\includegraphics[width=0.49\textwidth]{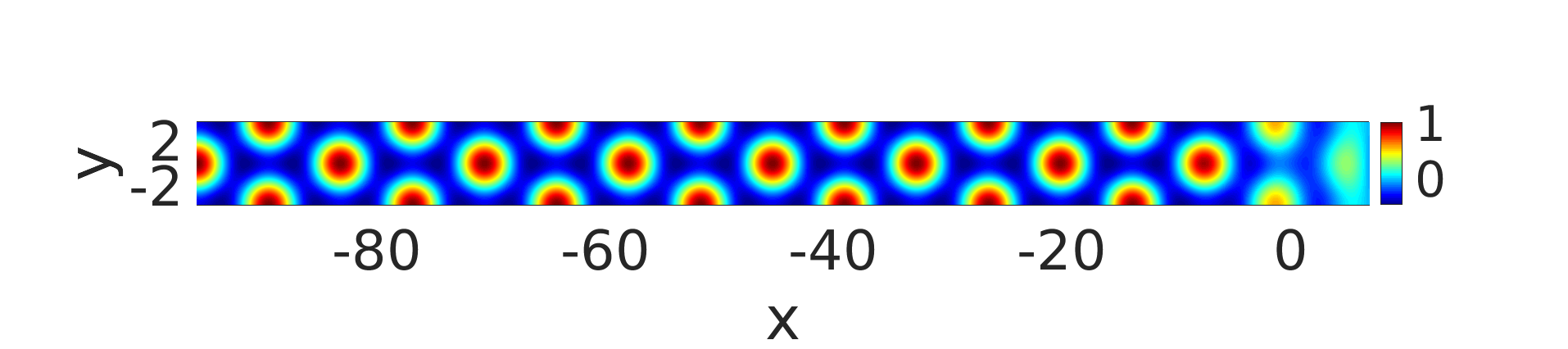}
\includegraphics[width=0.49\textwidth]{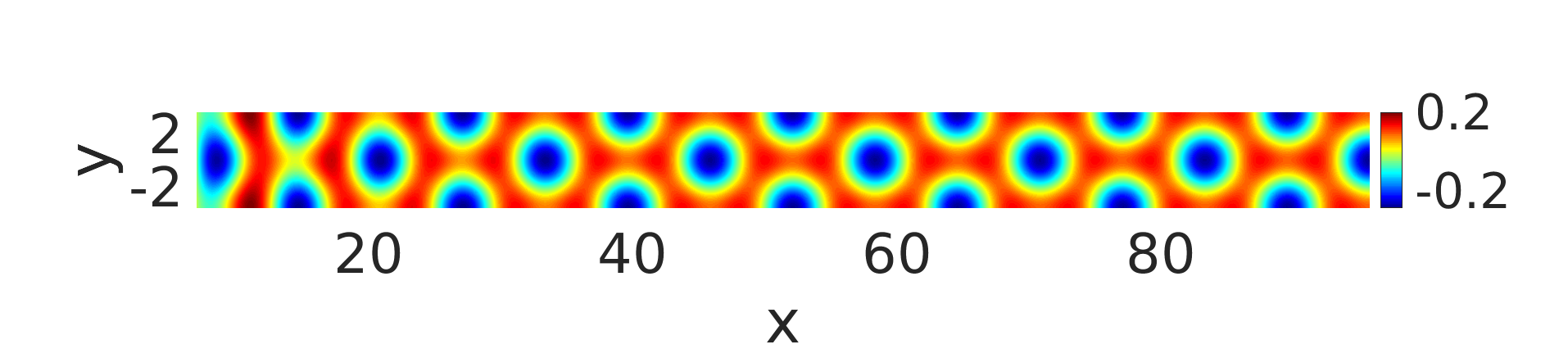}
\end{minipage}

\caption{(a) (a) Numerically calculated solution branches of \eqref{sh} for $\ep=0.5$ and $c_2=-2$, $c_3=-1$, $c_4=5$, $c_5=-2$ (as for \figrefb{stabregion} and \figref{hexsnake}). Here $\Omega=(-l_x,l_x)\times (-l_y,l_y)$ with $l_x=30.5\pi$, $l_y=2\pi / \sqrt{3}$ and we use Neumann boundary conditions. (b)-(i) Density plots of solutions labeled on the gray branch. (j) Left and right part of solution 458 plotted seperately. } \label{hexsnakeunsym}
\end{figure}

We used $n=61$ (odd) for developing \figref{hexsnakeunsym} and Neumann boundary conditions. In this case stripes corresponding to the wavenumber $k=1$ have a maximum and minimum on the left and right boundary, respectively, or vice versa. Hexagons and beans corresponding to the wavenumber $k=1$ cannot fulfill the boundary conditions, but a state, which passes near up- and down-hexagons such as solution 235 of \figref{hexsnake}, does. When the stripe branch loses its stability, it does not bifurcate a bean branch, but a branch like the gray one of \figrefa{hexsnake}.    

We describe the behavior of this branch in the following. 
A formation of a spot occurs on the left side (from bifurcation to solution 53). Then the branch turns around in a fold in the point, where stripes become stable. It follows a part where a gap front invades the stripes along the branch. Solution 80 is a state which is a mix of stripes, up-, and down-hexagons. Solution 120 consists of up- and down-hexagons and an interface which lies on the left side. It follows a part of the branch where it moves back and forth in parameter space by passing stable and unstable ranges (snaking). Along this snake the interface moves to the right side (the up-hexagons invade the down-hexagons). The interface shift is of the length $4 \pi$ from a fold on one side to the fold after next on the same side (compare solution 410 and 458). Because of the boundary problems described above the branch cannot go back to the spot branch as the gray branch of \figrefa{hexsnake} does, but terminates in a branch which is nearby (see solution 810). In \figrefj{hexsnakeunsym} we can see that we have no boundary problems for connections between up- and down-hexagons as we had in \figref{hexsnake}. Furthermore, one should mention that the gray branches of \figref{hexsnake} and \figref{hexsnakeunsym} show a similar behavior. The length of the snake in \figref{hexsnakeunsym} is larger than the one of \figref{hexsnake}, which is because we used a larger domain here.   

\subsection{Explanations via conserved quantities}

In the following we explain the behavior of the gray branches of \figref{hexsnake} and \figref{hexsnakeunsym} analytically via amplitude equations and Ginzburg-Landau energy techniques (see \cite{hexsnake,schnaki}) and numerically via the full system and the Hamiltonian (see \cite{hexsnake}). We start with the first way. The amplitude reduction \eqref{gl5} is a Landau system of the form 
\alinon{
\partial_t A_1&=f_1(A_1,A_2,A_3),\\
\partial_t A_2&=f_2(A_1,A_2,A_3),\\
\partial_t A_3&=f_3(A_1,A_2,A_3),
}
where $A_j=A_j(t)$ is time dependent. Assuming that $A_j=A_j(t,x)$ is also space dependent we obtain the Ginzburg-Landau system 
\ali{
\partial_t A_1&=c_0\partial_x^2A_1+f_1(A_1,A_2,A_3),\nonumber\\
\partial_t A_2&=\frac{c_0}{4}\partial_x^2A_2+f_2(A_1,A_2,A_3),\label{gls}\\
\partial_t A_3&=\frac{c_0}{4}\partial_x^2A_3+f_3(A_1,A_2,A_3).\nonumber
}
Background on this formal procedure and for so called attractivity and approximation theorems can be found in \cite{gs94,bsvh95,gs99b,alex,SU17}. 

Restricting to real amplitudes $A_1,  A_2,  A_3$,  we can split the total energy into the kinetic and potential energy, i.e., $E_{\text{t}}=E_{\text{k}}+E_{\text{p}}$, where
\ali{
E_{\text{k}}= & \ \frac{c_0}{2}\left( (\partial_xA_1)^2+\frac14(\partial_xA_2)^2+\frac14(\partial_xA_3)^2    \right),\nonumber \\ 
E_{\text{p}}= & \ \frac12 c_1 (A_1^2+A_2^2+A_3^2)+2c_2A_1A_2A_3\nonumber\\
&+\frac34 c_3 (A_1^4+A_2^4+A_3^4)+3 c_3 (A_1^2A_2^2+A_1^2A_3^2\nonumber\\
&+A_2^2A_3^2)+12c_4G_4+c_5G_5 \label{Ep}
}
with 
\alinon{
G_4=A_1^3A_2A_3+A_1A_2^3A_3+A_1A_2A_3^3
}
and
\alinon{
G_5=& \ 75A_1^2A_2^2A_3^2+\frac53(A_1^6+A_2^6+A_3^6)\\
& +15(A_1^2A_2^4+A_1^4A_2^2+A_1^2A_3^4+A_1^4A_3^2\\
&+A_2^2A_3^4+A_2^4A_3^2).
}

The Hamiltonian for stationary states of generalized Swift-Hohenberg equations over 2D domains can be found in \cite{hexsnake}. For \eqref{sh} this Hamiltonian is given by
\begin{equation}\label{ham}
\begin{aligned}
H(u)= & \int_{-l_x}^{l_x} \bigg[\frac{u_{xx}^2}{2}-u_{xxx}u_x-u_{x}^2+u_{xy}^2\\
&+u_{y}^2-\frac{u_{yy}^2}{2}+F \bigg] \text{d}y
\end{aligned}
\end{equation}
with
\ali{F=\frac{c_1-1}{2}u^2+\frac{c_2}{3}u^3+\frac{c_3}{4}u^4+\frac{c_4}{5}u^5+\frac{c_5}{6}u^6.
}

It holds $\frac{\text{d}}{\text{d}x}E_\text{t}=0$ and $\frac{\text{d}}{\text{d}x}H=0$, i.e., $E_\text{t}$  and $H$ are conserved. Thus, a necessary conditions for a heteroclinic front connection between two stationary states in \eqref{gls} and \eqref{sh} is that both have the same potential energy $E_\text{p}$ and Hamiltonian $H$, respectively. Using the Ginzburg-Landau energy is an easy and rapid way, but is only valid near onset. Using the Hamiltonian is more costly, but gives accurate results not only near onset.

\begin{figure}[h]
\begin{minipage}{0.23\textwidth}
(a)
\includegraphics[width=1\textwidth]{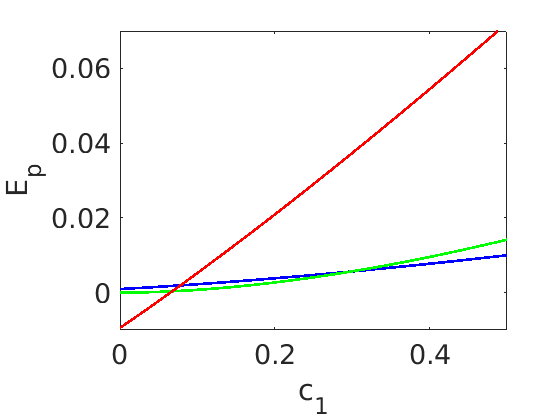}
\end{minipage}
\begin{minipage}{0.23\textwidth}
(b)
\includegraphics[width=1\textwidth]{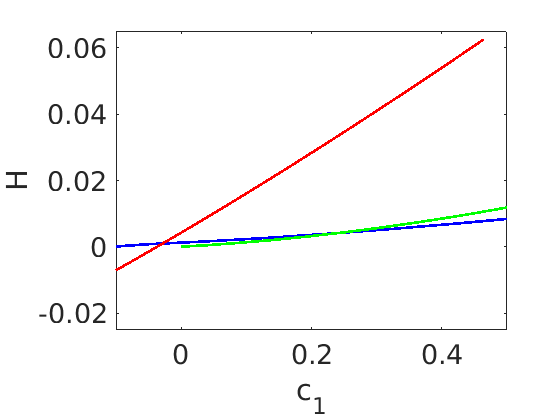}
\end{minipage}
\caption{Potential energy $E_\text{p}$ \eqref{Ep} and Hamiltonian for up-hexagons (red), stripes (green), and down-hexagons (blue) are shwon in (a) and (b) respectively.}\label{Epot}
\end{figure} 

The potential energy \eqref{Ep} and the Hamiltonian \eqref{ham} for stripes, up-, and down-hexagons are shown in \figrefa{Epot} and (b), respectively. First of all we can see that up- and down-hexagons have the same potential energy $E_\text{p}$ and Hamitonian $H$ for $c_1\approx 0$, which explains the existence of connections between up- and down-hexagons in this range. Furthermore, we can conclude that the state 53 on the gray branch of \figref{hexsnakeunsym} is not a state for which the corresponding orbit passes near stripes and up-hexagons, but we expect that the orbit moves from stripes into the spot direction and is never close to up-hexagons before moving back to stripes. 

Stripes and up-hexagons have the same potential energy $E_\text{p}$ and Hamitonian $H$ for $c_1\approx 0.27$ such that a necessary condition for the existence for a connection between stripes and up-hexagons is fulfilled. We did not try to find such a state for studying the snaking behavior of the branch, since this is already done in \cite{schnaki,w16} for other systems and we do not expect to find a different behavior. But what we can see is that the invasion of down-hexagons between solution 80 and 120 starts when the potential energy of stripes becomes smaller than the one of down-hexagons.  

\section{Conclusion}
We cannot predict bistable ranges between up- and down-hexagons and tristable ranges between stripes, up- and down-hexagons if we use a third order amplitude reduction for pattern forming systems. If we use fifth order amplitude reductions, we are able to find such ranges. We studied a generalized Swift-Hohenberg equation \eqref{sh} which is scaled with a parameter $\ep$, which makes a fifth order amplitude reduction easy and valid for small $\ep$. We fixed $c_3,c_4,c_5$ for finding tristable ranges between stripes, up- and down-hexagons in the $c_1$-$c_2$-plane based on the amplitude reduction. We fixed $c_2$ and used $c_1$ as bifurcation parameter for calculating branches of stripes, up- and down-hexagons for $\ep=0.5$ and $\ep=1$ using numerical path following methods. We compared these branches with the one obtained from the amplitude reductions and saw that we obtain acceptable predictions for $\ep=0.5$ and a completely different bifurcation diagram for $\ep=1$. Here we found ranges, where five different states are stable.

We went on with $\ep=0.5$ for finding interesting states existing in the multistable ranges and found a snaking branch of connections between up- and down-hexagons. This branch bifurcates from the stripe branch and terminates on the spot branch. We used conserved quantities in order to understand the location of the snaking part.

\bibliography{snbib}{}
\bibliographystyle{plain}

\end{document}